\documentclass[12pt,leqno]{article} 
\usepackage{amsmath,amsfonts,mathabx}

\newtheorem{theorem}{Theorem}

\newtheorem{definition}{Definition}
\newtheorem{corollary}[theorem]{Corollary}

\begin{document}

\title{\textbf{Sasakian structures \\ A foliated approach}}

\author {Robert A. Wolak \\{\small Wydzial Mathematyki i Informatyki, Uniwersytet Jagiellonski}\\{\small Lojasiewicza 6,  30-348 Krakow, Poland, robert.wolak@im.uj.edu.pl.}\\
\\ Dedicated to Micha\l \;   Heller on his 80th birthday
}

\maketitle

Recent renewed interest in Sasakian manifolds is due mainly to the fact that they can provide examples of generalized Einstein manifolds, manifolds which are of great interest in mathematical models of various aspects of physical phenomena.  Sasakian manifolds are odd dimensional counterparts of K\"ahlerian manifolds to which they are closely related. The book of Ch. Boyer and K. Galicki, {\it Sasakian Geometry \/} is both the best introduction to the subject and at the same time it  gathers state of the art information and results on these manifolds. However, although the authors are well aware that a Sasakian structure is a very special one-dimensional Riemannian foliation with K\"ahlerian transverse structure, they use this fact only in a few very special cases. 

The paper  presents an approach to Sasakian manifolds on which the author gave several lectures, most recently at the  {\it Workshop on almost hermitian and contact geometry\/} at the Banach Center in B\c{e}dlewo in October 2015 and at University of the Basque Country in February 2016.  The first lectures on the topic the author gave  at Universidad de Sevilla in October 1988 and then presented the consequence for the geometry of Sasakian manifolds, in particular the relationns between various curvatures and those of the transverse K\"ahler manifold. The results were published in several sections of  \cite{WO_S} as well as in \cite{Wo_deb}. The most general theory of geometrical structures "adapted" to a foliation was presented in \cite{WO_T}, see also \cite{WO_S}. The paper concentrates on cohomological properties of Sasakian manifolds and of  transversely holomorphic and K\"ahlerian foliations. These properties permit to formulate obstructions to the existence of Sasakian structures on compact manifolds. The presented results are due to the author as well as his former and present Ph.D. students.

\section{Geometric structures on contact manifolds}

In this paper we are going to work on odd dimensional smooth manifolds. Let $M$ be a smooth connected manifold of dimension $m=2n+1.$

Let $\eta$ be a 1-form on $M$ such that $\eta \wedge d\eta^n \neq 0,$ i.e., it is a volume form of $M.$ Then $\eta$ is called a {\it contact form.} Two contact forms $\eta$ and $\eta'$ are said to be equivalent if there exists a smooth function $f$ such that $\eta = f\eta' .$ Such an equivalence class is called a contact structure. The pair $(M, \eta)$ usually is called a {\it strict contact manifold.} A contact manifold is a smooth manifold $M$ with a contact structure $[\eta ].$

\medskip

On a strict contact manifold $(M, \eta)$ there exists a unique vector field $\xi ,$ called the {\t Reeb vector field,} such that

$$ \eta (\xi ) = 1 \;\;\; \wedge \;\;\; i_{\xi} d\eta = 0.$$

\noindent
This property ensures that the 1-dimensional foliation ${\mathcal F}_{\xi}$ generated by the non-vanishing vector field $\xi$ is transversely symplectic. This foliation is independent of the choice of the contact form $\eta$ within the equivalence class, it is one of the objects we can associate to  a contact manfold. On the contact manifold $(M,  [\eta ])$ we have the canonical splitting of the tangent bundle $TM :$

$$ TM = T{\mathcal F}_{\xi} \oplus D$$

\noindent
where $D = ker \eta .$

\medskip
Geometers are more at home with a richer structure: {\it almost contact (manifold)\/}:

\begin{definition}
An almost contact structure on a smooth manifold M is a triple $(\xi , \eta, \phi )$ where

i) $\xi $ is a vector field on M,

ii) $\eta$ is a 1-form on M,

iii) $\phi$ is an endomorphism of the tangent bundle TM

\noindent
such that

$$ \eta (\xi) = 1, \;\;\;   \;\;\; \phi^2 = - id_{TM} + \xi \otimes \eta .$$

\end{definition}

\noindent
One can easily verify that $\phi (\xi ) = 0,$ $\eta \phi = 0,$ and that the  tangent bundle $TM$ splits naturally into the direct sum ${\mathcal F}_{\xi} \oplus D$ where $D = ker\eta = im \phi .$

\medskip
The next step in the enrichment of the geometrical structure is a compatible Riemannian metric.

\begin{definition}

A Riemannian metric g is said to be compatible with an almost contact structure $(\xi , \eta, \phi )$ if for any vector fields X and Y

$$ g( \phi (X), \phi (Y) ) = g(X,Y) - \eta (X) \eta (Y). $$

\noindent
Then the quadruple $(g, \xi , \eta, \phi )$ is called an almost contact metric structure.

\end{definition}

\noindent
{\bf Remark} The vector field $\xi$ need not be a Killing vector field for the metric $g$ even if the  quadruple $(g, \xi , \eta, \phi )$ is  an almost contact metric structure, cf. \cite{BG}.

\medskip

Combining the topological and geometrical structures we get the so-called {\it contact metric structure\/}:

\begin{definition}

An almost contact structure $(\xi , \eta ', \phi )$ is said to be compatible with the strict contact structure (form)  $\eta$ if $\eta = \eta ',$ $\xi$ is its Reeb vector field, and for any vector fields $X$ and $Y$

$$ d\eta (\phi (X), \phi (Y)) = d\eta (X,Y) \;\;\; and \;\;\; d\eta (\phi (X), X) > 0, X \in D, \;\;\;  X \neq  0 .$$

\end{definition}

\begin{definition}

A strict contact manifold $(M, \eta)$  with a compatible  almost contact metric structure $(g, \xi , \eta ', \phi )$ such that for any two vector fields X and Y

$$g(X, \phi (Y)) = d\eta (X,Y) $$

\noindent
is called a contact metric structure.

\end{definition}

The Reeb vector field of a contact metric structure need not be Killing.  If it is,  the structure is called K-contact.

\begin{definition}

A  contact metric structure  $(g, \xi , \eta ', \phi )$  on the manifold $M$ is called  K-contact if its Reeb vector field $\xi$ is a Killing vector field of the Riemannian metric $g.$ Then  $(M, g, \xi , \eta ', \phi )$ is called a K-contact manifold.
\end{definition}

And finally, the most complex structure considered is the Sasakian structure (manifold).

\begin{definition}
An almost contact structure $(\xi , \eta ', \phi )$ on the manifold M is normal iff

$$N_{\phi}(X,Y) = [\phi X, \phi Y] + \phi^2[X,Y] - \phi [X, \phi Y] - \phi [\phi X, Y] = -2\xi \otimes d\eta (X,Y)$$
\noindent
for any  vector field X and Y on M.
\end{definition}

\begin{definition} 
A K-contact manifold whose underlying almost contact structure is normal is called a Sasakian manifold.
\end{definition}


\section{Transverse properties of Sasakian manifolds}

Let $\mathcal F$ be a foliation on a Riemannian $m$-manifold $(M,\, g).$ Then $\mathcal F$ is defined by a cocycle    $\mathcal U =\{ U_i , f_i , g_{ij} \}_{i,j\in I}$ modeled on a $2q$-manifold $N_0$ such that 

(1) $\{ U_i \}_{i\in I}$ is an open covering of M, 

(2) $f_i : U_i \to N_0$ are submersions with connected fibres, 

(3) $g_{ij} : N_0\to N_0 $ are local diffeomorphisms of $N_0$ with $f_i = g_{ij} f_j$ on $U_i \cap U_j .$ 

\noindent
The connected components of the trace of any leaf of $\mathcal F$ on $U_i$ consist of the fibres of $f_i .$ The open subsets
$N_i = f_i ( U_i )\subset N_0$ form a $q$-manifold $N_{\mathcal U}=\amalg N_i$,
which can be considered as a transverse manifold of the
foliation $\mathcal F .$ The pseudogroup $\mathcal H_{\mathcal U}$ of local
diffeomorphisms of $N$ generated by $g_{ij}$ is called the holonomy
pseudogroup of the foliated manifold $(M, \mathcal F )$ defined by the
cocycle $\mathcal U .$ Different cocycles can define the same foliation, then we have two different transverse manifolds and two holonomy pseudogroups. In fact,  these two holonomy  pseudogroups are equivalent in the sense of Haefliger, cf. \cite{Hae}. 

According to Haefliger, cf. \cite{Hae_80}, a transverse property of a foliated manifold is a property of foliations which which is shared by any two foliations with equivalent holonomy pseudogroup. For example, being Riemannian, transversely symplectic, transversely almost-complex, transversely K\"ahler, etc., is a transverse property. A Riemannian foliation, i.e., admitting a bundle-like metric, is defined by a cocycle $\mathcal U$ modelled on a Riemannian manifold whose local submersions are Riemannian submersions. Then the associated transverse manifold $N_{\mathcal U}$ is Riemannian and the associated holonomy pseudogroup ${\mathcal H}_{\mathcal U}$ is a pseudogroup of local isometries. Any foliation defined by a cocycle $\mathcal V$ whose holonomy psudogroup ${\mathcal H}_{\mathcal V}$ is equivanet to ${\mathcal H}_{\mathcal U}$ is also Riemannian, as the equivalence of pseudogroups transports the Riemannian metric from $N_{\mathcal U}$ to $N_{\mathcal V}$ and ensures that the pseudogroup ${\mathcal H}_{\mathcal V}$ is a pseudogroup of local isometries of the transported metric. This metric can be lifted to a bundle-like metric (not unique) on the other foliated manifold making the second foliation Riemannian. The same procedure can be applied to any geometrical structure, for the discussion of this general procedure see \cite{WO_T,WO_S}.

\medskip

The space of  ${\mathcal H}_{\mathcal U}$-invariant $k$-forms on the manifold  $N_{\mathcal U}$ can be identified with the space of foliated sections of the bundle $\bigwedge ^kN(M,{\mathcal F})^*$  which in turn is isomorphic to the space of $k$ basic forms 

$$A^k(M,{\mathcal F}) = \{ \alpha \in A^k(M) \colon i_X\alpha = i_Xd\alpha =0 \;\; for \;\; all \;\; vectors\;\;  X  \in {\mathcal F}\}$$

\noindent
The differential sends basic forms to basic forms and the cohomology of the complex $(A^*(M,{\mathcal F}), d)$ is called the {\it basic cohomology\/}  of the foliated manifold $(M,{\mathcal F}).$ In the language of basic cohomology we can express a very important property of foliations:

\begin{definition} A foliation $\mathcal{F}$ on M is called homologically orientable if $H^{cod\mathcal{F}}(M,\mathcal{F})=\mathbb{R}$.
\end{definition}

For the discussion the meaning and importance of the condition see \cite{El_2010}

\medskip

   Let  $\phi : U\to R^p \times R^q ,\,\phi =(\phi^1 ,\phi^2)=(x_1, ..., x_p ,y_1,...,y_q )$ be an adapted chart on a
   foliated manifold $(M, \mathcal F)$. Then on $U$   the  vector fields
   $\frac{\partial}{\partial x_1},...\frac{\partial}{\partial x_p}$ span
   the bundle $T\mathcal F$ tangent to the leaves of the foliation $\mathcal F ,$ the equivalence
   classes denoted by $\bar{\frac{\partial}{\partial y_1}}, ...\bar
   {\frac{\partial}{\partial y_q}}$ of  $\frac{\partial}{\partial y_1}, ...
   \frac{\partial}{\partial y_q}$ span the normal bundle $N(M, \mathcal F )
   =TM / T\mathcal F $ which is isomorphic to the subbundle $T\mathcal F^\perp$. 

\medskip

All the definitions of Section 1 have been formulated in a purely geometrical way without any reference to the characteristic foliation. Let us look at the transverse structure of the characteristic foliation.

The characteristic foliations of a contact manifold $(M, [ \eta ])$ is transversely symplectic as the 2-form $d\eta $ is basic and defines a transverse symplectic form. 

The basic cohomology class $[ d\eta ] \in H^2 (M, {\mathcal F})$ is in the kernel of the natural mapping $H^2 (M, {\mathcal F}) \rightarrow H^2(M), $ and $[ d\eta ]^n \in H^{2n} (M, {\mathcal F})$ is in the kernel of the natural mapping $H^{2n} (M, {\mathcal F}) \rightarrow H^{2n}(M). $
Therefore if transverse volume form  $[\omega ]^n$ defines a non-zero basic cohomology class, then this 2n-form is in the kernel  the natural mapping $H^{2n} (M, {\mathcal F}) \rightarrow H^{2n}(M), $ thus this mapping cannot be injective providing an obstruction to a transversely symplectic 1-dimensional foliation being the characteristic foliation of a contact structure. 

For example, such a 1-dimensional foliation cannot admit a transverse foliation with a compact leaf, as then according to the result of Molino-Sergiescu this mapping should be injective, cf. Theorem 2 of \cite{MS}.

\medskip


In the case of an almost contact structure $(\xi , \eta, \phi )$ on a smooth manifold $M,$ the following conditions are equivalent, cf. \cite{G} or Lemma 6.3.3 of \cite{BG}:

1) there exists a Riemannian metric for which the orbits of $\xi$ are geodesics;

2) $L_{\xi}\eta =0:$

3) $i_{\xi}d\eta =0.$

\medskip
\noindent
The conditions (2) and (3) are evidently equivalent, and they just say that the 2-form $d\eta$ is basic. 

\medskip


In the case of a contact metric structure  $(g, \xi , \eta ', \phi )$ we have the following equivalent conditions (Proposition 6.4.8 of \cite{BG})

1) the characteristic foliation is Riemannian for $g;$

2) the metric $g$ is bundle-like;

3) the vector field $\xi$ is Killing;

4)  the vector field $\xi$ preserves the (1,1)-tensor field $\phi ,$ i.e., $L_{\xi} \phi =0;$

5)  the  contact metric structure  $(g, \xi , \eta ', \phi )$ is K-contact.


\medskip
\noindent
Therefore   the characteristic foliation of a K-contact manifold   is transversely almost-K\"ahler, and the characteristic foliation of a Sasakian manifold is transversely K\"ahler, \cite{BG} Theorem 7.1.3. However, even if the characteristic foliation of a K-contact manifold is transversely K\"ahler, it does not imply that the structure is Sasakian. 

\section{Transversely K\"ahler foliations}

We have noticed that the characteristic foliation of a Sasakian manifold is transversely K\"ahler. In this section we will gather the results which are particular to transversely K\"ahler foliations and therefore are also true for the characteristic foliations of a Sasakian manifold. It will facilitate the search for characterizations of Sasakian manifolds and properties which can distinguish between K-contact and Sasakian manifolds.  

Let ${\mathcal F}$ be a foliation of dimension $p$ and codimension $2q$ on a smooth manifold $M$ of dimension $m=p+2q.$ It is a transversely K\"ahler foliation if there is a cocycle ${\mathcal U} = \{ (U_i,f_i,g_{ij})\}_{i,j \in I}$ defining the foliation $\mathcal F$ modelled on a K\"ahler manifold $(N,g_N, J_N)$ such that the local diffeomorphisms $g_{ij}$  of $N$ are K\"ahler isometries, or equivalently that  the associated holonomy pseudogroup ${\mathcal H}_{\mathcal U}$ is a pseudogroup of K\"ahler isometries of a K\"ahler structure on the transverse manifold 
$N_{\mathcal U}.$ 

We assume that the foliation $\mathcal F$ is transversely holomorphic, of complex codimension $q$ and that the manifold $M$ is compact.
 Therefore on the normal bundle $N(M,{\mathcal F})$ of the foliation ${\mathcal F}$ we have a foliated K\"ahler structure, i.e. a foliated Riemannian metric $\bar{g}$ and an endomorphism $\bar{J}$ of the normal bundle such that $\bar{J}^2 =-Id$ (an almost complex structure in the normal bundle) compatible with $\bar{g},$ i.e. for any $X,Y \in N(M,{\mathcal F})$

$$\bar{g}(\bar{J}X, \bar{J}Y) = \bar{g}(X,Y)$$

\noindent
and satisfying the " integrability " condition:  for any $X$ and $Y$

$$N_{\bar{J}}(X,Y) = [\bar{J}X,\bar{J}Y] +\bar{J}^2[X,Y] -\bar{J}[\bar{J}X,Y] -\bar{J}[X,\bar{J}Y] =0$$

\noindent
Then the fundamental  K\"ahler 2-form $\bar{\Omega}(X,Y) = \bar{g}(X,\bar{J}Y)$ is basic and correspond to the K\"ahler form of the transverse K\"ahler structure of the foliation. For any point $x$ of $M$  there exists an adapted chart $(x_1,...,x_{m-2q}, z_1,...,z_q)$ modelled on $R^{m-2q} \times C^q$ defined on an open neighbourhood of $x.$ Basic forms on $(M,{\mathcal F})$ are in one-to-one correspondence with holonomy invariant ($\mathcal H$-invariant) forms on the transverse manifold. Basic $k$-forms are just foliated sections of the kth exterior product of the conormal bunsle $N(M,{\mathcal F})^*$ foliated by the natural lift of the foliation $\mathcal F$, cf. \cite{WO_T}. If the foliation $\mathcal F$ is transversely holomorphic its normal bundle has a (almost) complex structure corresponding to the complex structure of the transverse manifold. Therefore any complex valued basic $k$-form can be represented as a sum of $k$-forms of pure type $(r,s)$  corresponding to the decomposition of $k$-forms on the complex manifold $N.$ It is possible as elements of the holonomy pseudogroup are bi-holomorphic local diffeomorphisms. This decomposition can be also obtained by considering the decomposition of sections of the complex bundle $\Lambda^k_{C}N(M,{\mathcal F})^*.$ Then a basic $k$-form $\alpha$  is of pure type $(r,s)$ if for any point of $M$ there exists an adapted chart  $(x_1,...,x_{m-2q}, z_1,...,z_q)$ such that

$$\alpha  = \Sigma f_{IJ}dz_{i_1} \wedge ... \wedge dz_{i_r}\wedge d\bar{z}_{j_1} \wedge ... \wedge d\bar{z}_{j_s}$$
 
\noindent
where  $1 \leq i_1 < ...  < i_r \leq q, 1 \leq j_1 < ...  < j_s \leq q, I=(i_1,...,i_r), J = (j_1, ..., j_s).$

\medskip

Let us denote by $A^k_C(M,{\mathcal F})$ the space of complex valued basic $k$-forms on the foliated manifold $(M,{\mathcal F}),$ and by $A^{r,s}_C(M,{\mathcal F})$ the space of complex valued basic forms of pure type $(r,s).$ Then

$$ A^k_C(M,{\mathcal F})  =  \Sigma_{r+s=k} A^{r,s}_C(M,{\mathcal F})$$

\noindent
for short

$$A^k =  \Sigma_{r+s=k} A^{r,s}.$$

\medskip
\noindent
The exterior differential $ d \colon A^k_C(M,{\mathcal F}) \rightarrow A^{k+1}_C(M,{\mathcal F})$ is decomposed into two components 
$\partial$ and $\bar{\partial}$ of bidegree $(1,0)$ and $(0,1),$ correspondingly,

$$  \partial \colon A^{r,s} \rightarrow A^{r+1,s} \;\;\; and \;\;\;  \bar{\partial} \colon A^{r,s} \rightarrow A^{r,s+1}.$$

\medskip

The basic cohomology of transversely holomorphic and transversely K\"ahler foliations was studied in depth by A. El Kacimi-Alaoui, cf. \cite{E1}. We recall some basic results from this paper. 

\medskip

We assume that the foliation is transversely Hermitian. The operator

$$ * \colon A^k(M,{\mathcal F}) \rightarrow A^{2q-k}(M,{\mathcal F})$$

\noindent
defined in \cite{ElH} using the transverse part of the bundle-like metric, and the corresponding  standard * operator on the level of the transverse manifold, can be extended to an operator

$$ \bar{ *} \colon A^k_C (M,{\mathcal F}) \rightarrow A^{2q-k}_C (M,{\mathcal F}).$$

\medskip
\noindent
The normal part of the bundle-like metric $g$, or the corresponding transverse metric, defines a Riemannian (Hermitian) metric $g^k$  on the bundle $\Lambda^k_{C}N(M,{\mathcal F})^*,$ and therefore we can define a scalar product on 

$$ A^*_C(M,{\mathcal F}) = \Sigma A^k_C(M,{\mathcal F})$$

\noindent
as follows, \cite{El_2010}:

\medskip
\noindent
if $ \alpha \in A^k_C(M,{\mathcal F}),$   $\beta \in A^l_C(M,{\mathcal F}),$  and $k \neq l$

$$<\alpha , \beta > = 0$$

\medskip
\noindent
if $ \alpha \in A^k_C(M,{\mathcal F}),$   and $\beta \in A^k_C(M,{\mathcal F}),$ 

$$<\alpha , \beta > =  \int_M g^k ( \alpha, \beta )$$

\noindent
The operator $\delta \colon  A^k_C(M,{\mathcal F}) \rightarrow A^{k-1}(M,{\mathcal F})$ defined as $ \delta = \bar{*}d\bar{*}$ is the adjoint operator of $d$ with respect to the scalar product $<,>$. 

\medskip

Following the classical (manifold) case we define the "foliated" Laplacian 

$$ \Delta  = d\delta + \delta d.$$

\noindent
The foliated Laplacian sends basic forms into basic forms, it is a self-adjoint foliated (transversely) elliptic operator, cf. \cite{E1}.

\medskip

We can also define basic Dolbeault cohomology of the foliated manifold $(M,{\mathcal F}).$ For a fixed $r, $ $0 \leq r \leq q,$ consider the differential complex:

$$ 0 \rightarrow A^{r,0}  \xrightarrow{\bar{\partial}}  A^{r,1} \xrightarrow {\bar{\partial}} ...  
 \xrightarrow {\bar{\partial}}  A^{r,q} \xrightarrow {\bar{\partial}} 0.$$

\noindent
Its cohomology is called the  basic Dolbeault cohomology of the foliated manifold $(M,{\mathcal F}),$ and denoted

$$ H^{r,*} (M,{\mathcal F})  = \Sigma_{s=0}^q H^{r,s} (M,{\mathcal F})  .$$

\noindent
The operator $\bar{*}$ induces  an  isomorphism  $\bar{*} \colon A^{r,s} \rightarrow A^{q-r,q-s}.$ Using the same procedure as for the operator $\delta$ we define an operator $\bar{\delta}$ by the formula

$$ \bar{\delta} = - \bar{*} \bar{\partial} \bar{*}.$$

\noindent
The operator $ \bar{\delta}$ is the adjoint of $\bar{\partial}$ with respect to the just defined scalar product. Moreover, the operator

$$\Delta" = \bar{\partial}\bar{\delta}  + \bar{\delta}\bar{\partial}$$

\noindent
is a self-adjoint foliated (transversely) elliptic operator.

\medskip

In the case of transversely K\"ahler foliations we can say much more about the basic cohomology and operators just defined.

\medskip

The K\"ahler form of the transverse manifold $N$ corresponds to a basic $(1,1)$-form on $(M,{\mathcal F})$ which we call the (transverse) K\"ahler form of the foliated manifold. Using this form we define the $L$ operator

$$ L \colon A^k_C(M,{\mathcal F}) \rightarrow A^{k+2}_C(M,{\mathcal F})$$

$$ L\alpha = \alpha \wedge \omega.$$

\noindent
Its adjoint with respect to $<,>$ is $\Lambda  = - \bar{*} L \bar{*}.$

\medskip

For transversely K\"ahler foliations on compact manifolds we have the following relations:

$$\Lambda \partial - \partial \Lambda = - \sqrt{-1} \bar{\delta},$$

$$\Lambda \bar{ \partial}  - \bar{\partial} \Lambda = - \sqrt{-1} {\delta},$$

$$ \partial \bar{\delta} + \bar{\delta}\partial =  \delta \bar{\partial} + \bar{\partial} \delta =0,$$

$$ \Delta = 2\Delta ",$$

$$\Delta L = L \Delta, \;\;\; \Delta \Lambda = \Lambda \Delta.$$

 These identities permited A. ElKacimi Alaoui to prove the following theorem, cf. \cite{E1}

\begin{theorem} Let $\mathcal F$ be a transversely K\"ahler foliation on a compact manifold M. If $\mathcal F$ is homologically oriented, then

i) a basic k-form $\alpha = \Sigma_{r+s=k}\alpha_{r,s}, \;\; \alpha_{r,s} \in A^{r,s}$ is harmonic if and only if the forms $\alpha_{r,s}$ are harmonic, thus 

$$ H^k_C (M,{\mathcal F}) \cong  \Sigma_{r+s=k} H^{r,s}(M,{\mathcal F}).$$

ii) the conjugation induces isomorphisms $$ H^{r,s}(M,{\mathcal F}) \cong H^{s,r}(M,{\mathcal F}).$$

iii) for any $0 \leq r \leq q,$  the form $\omega^r$ is harmonic, thus $H^{r,r}(M,{\mathcal F}) \neq 0.$

\end{theorem}

The complex $ A^* = \Sigma_{r,s} A^{r,s} $ of complex valued basic forms can be filtered by

$$F^kA =  \Sigma_{r \geq k} A^{r,*}.$$

\noindent
The filtration is compatible with the bigradation of the complex. Therefore we can define the associated spectral sequence which is called the basic Fr\"olicher spectral sequence of the transversely holomorphic foliation $\mathcal F$ , cf. \cite{CW_Pa} . It converges to the complex basic cohomology of the foliated manifold $(M,{\mathcal F}).$ 

\medskip

The terms $E^{r,s}_1$ are just the basic $(r,s)$-Dolbeault  cohomology groups. If the foliation $\mathcal F$ is homologically oriented and transversely K\"ahler,  then it is a simple consequence of Theorem 1 that the Dolbeault spectral sequence collapses at the first term, cf. Theorem 2 of \cite{CW_Pa}. Indeed, the Hodge theorem combined with the just metioned theorem ensures that

$$ E^{r,s}_1 \cong {\mathcal H}^{r,s} (M,{\mathcal F})$$

\noindent
where ${\mathcal H}^{r,s} (M,{\mathcal F})$  is the space of $(r,s)$-pure harmonic forms. As harmonic forms are closed the operator $d_1$ is trivial(vanishes). 

\medskip

In \cite{CW_Pa} the authors noticed that the so called $dd^c$- lemma   for K\"ahler manifolds is an algebraic consequence of several identities. These identities have their counterparts for the basic cohomology of the transversely K\"ahler foliation on a compact manifold so the  $dd^c$- lemma is also true for the basic cohomology of a transversely K\"ahler foliation on a compact manifold. On the other hand  the proof of the formality of the cohomology of a compact K\"ahler manifold, \cite{DGMS}. Therefore retracing the steps of the original proof we obtain

\begin{theorem}
Let $\mathcal F$  be a transversely K\"ahler foliation on a compact manifold M. If  $\mathcal F$ is homologically oriented then the minimal
model of the complex basic cohomology of $\mathcal F$ is formal and thus Massey products of complex valued basic forms vanish.
\end{theorem}

\section{Obstructions to existence of Sasakian structures}

We have remarked that the characteristic foliation of the Sasakian manifold is transversely K\"ahler. Therefore we have a 1-dimensional (tangentially) orientable foliation with a very sophisticated transverse structure.  Moreover, the normality condition is not a {\it transverse property} as its formulation involves vectors tangent to leaves of the foliation, in particular the characteristic vector field $\xi$. The corresponding transverse property can be formulated as follows:

Let $\bar{J} \colon N(M, {\mathcal F}) \rightarrow  N(M, {\mathcal F})$ be the endomorphism of the normal bundle defined for any tangent vector $X$ as

$$ \bar{J} (\bar{X} ) = \overline{ \phi (X)  }$$

\noindent
where $\bar{X}$ is the vector in the normal bundle corresponding to   a tangent vectoc $X.$ The  endomorphism is well defined and $ \bar{J}^2 = -id.$ Therefore it is an almost complex structure in the normal bundle. 

\medskip

The vector field $\xi$ acts on the normal bundle, and therefore on the endomorphism $\bar{J}.$ It is a foliated endomorphism iff 

$$L_{\xi}\bar{J} =0, $$

\noindent
i..e. iff $L_{\xi}\bar{J}(\bar{X}) = 0 = \overline{[\xi, \phi (X)]} - \overline{ \phi ([\xi , X])}$ for  any foliated vection of the distribution $D.$ Then $\bar{J}$ corresponds to an almost complex structure $J$ on the transverse manifold of the characteristic  foliation. The normality condition insures that  $J$ is integrable (i.e., the Nujenhuis tensor $N_J=0,$). However, the normality condition is stronger, the equality
$N_J=0$ equivalent to  $N_{\bar{J}}=0,$ and thus to the fact that for any sections  $X, Y $ of  $D,$ $N_{\phi} (X,Y) $ is a vector field tangent to the characteristic foliation ${\mathcal F}_{\xi},$ i.e., of the form $h \xi$ for some smooth function $h$ on $M,$ but not necessarily $2d\eta (X,Y )$ as requires the normality condition. 

\medskip

Therefore having given a 1-dimensional foliation we can ask many questions like: 

\medskip

{ \it Is this foliation Riemannian, (transversely) Hermitian, transversely symplectic, transversely holomorphic, transversely K\"ahler?} 

\medskip

These questions are about the transverse structure of the foliation and can be answered in the language of transverse properties, so the basic cohomology can provide some obstructions to the existence of such structures. 

\medskip

It is not difficult to see that a  1-dimensional transversely K\"ahler foliation admits a contact metric structure in the sense that  it is the characteristic foliation of this structure. If the manifold $M$ is compact, the non-triviality of the basic cohomology ensures that one can modify the Riemannian metric to ensure that the foliation is Riemannian and minimal, i.e. generated by a Kiling vector field.

\medskip

Let $\xi$ be a non-vanishing vector field on the manifold $M.$ Assume that the foliation ${\mathcal F}_{\xi}$ generated by $\xi$ is transversely K\"ahler. Therefore on the transverse manifold $N$ of the foliated manifold $(M, {\cal F}_{\xi})$ there exists a holonomy invariant K\"ahler structure $(\hat{g},\hat{J}),$ i.e. $\hat{g}$ is a Riemannian metric, $\hat{J}$ a complex structure, and for any tangent vectors $X.Y$ of $N$

$$\hat{g}(\hat{J}X,\hat{J}Y) = g(X,Y)  \;\;\; and \;\;\; \hat{ \Omega} (X,Y) =\hat{g}(X,\hat{J}Y) \;\;\; is \; a \; closed \; 2-form. $$

\noindent
We can lift the Riemannian metric $\hat{g}$ to a Riemannian metric $\bar{g}$ in the normal bundle by the formula

$$ \bar{g}_y (\bar{X}, \bar{Y}) = \hat{g}_{f_i(y)} ( df_i(\bar{X}), df_i(\bar{Y}))$$

\noindent
where $f_i \colon U_i \rightarrow N_0$ is a submersion from a cocycle defining the foliation ${\mathcal F}_{\xi},$ $\bar{X}, \bar{Y} \in N(M, {\mathcal F}_{\xi})_y,$ and $y\in U_i.$ The complex structure is lifted to an almost comlex structure $\bar{J}$ in the normal bundle in a similar way:

$$df_i(\bar{J}_y(\bar{X}) )= \hat{J}_{f_i(y)} (df_i(\bar{X})).$$

\noindent
Then the associated 2-form $\bar{\Omega}(\bar{X}, \bar{Y}) = \bar{g}(\bar{X}, \bar{J}\bar{Y}) $ is (locally) the pull-back of the K\"ahler 2-form $\hat{\Omega},$ i.e.

$$f^*_i\hat{\Omega} = \bar{\Omega}$$

\noindent
for any submersion $f_i$ from the cocycle defining the foliation ${\mathcal F}_{\xi}.$
Next choose a suplementary subbundle $D$ to the foliation, which is isomorphic to the normal bundle as a vector bundle,   and define the Riemannian metric $g$ on $M$ as follows:

\medskip
\noindent
the subbundles $T{\mathcal F}_{\xi}$ and $D$ are orthogonal, $g(\xi,\xi ) =1,$ and transport $\bar{g}$ via the isomorphism  to $D$.

\medskip
The tensor field $\phi$ is defined in a similar fashion:

\medskip

\noindent
for vectors from the subbundle $D$ we define $\phi$ as the pull-back of $\bar{J}$ via the isomorphism from the normal bundle, and $\phi (\xi) =0.$

\medskip
Let us define the 1-form $\eta$ as 

$$\eta (X) = g(\xi, X).$$

\noindent
Then, obviously,  the triple $(\xi ,\eta, \phi ) $ is an almost contact structure on the manifold $M.$  Let $X, Y$ be any vectors on $M.$ Taking into account the splitting $T{\mathcal F}_{\xi} \oplus D$ we can write $X = a_X\xi + \bar{X}$ and $Y = a_Y \xi + \bar{Y}.$ Thus

$$g(X,Y) = g(\bar{X},\bar{Y}) + g(a_X\xi , a_Y \xi) = g(\bar{X},\bar{Y})  + a_Xa_Y  = g(\bar{X},\bar{Y}) + \eta (X) \eta (Y) $$

$$\;\;  = g(\phi (\bar{X}), \phi (\bar{Y})) + \eta (X) \eta (Y) = g(\phi (X), \phi (Y)) + \eta (X) \eta (Y)  $$

\noindent
as the metric $\bar{g}$ is $\bar{J}$-invariant. Therefore the quadruple $(g, \xi , \eta , \phi )$ is an almost contact metric structure. 

\medskip

Assume that the vector field $\xi $ is the characteristic vector field of a (strict) contact structure $\eta$ on the manifold $M.$ The 2-form $d\eta$ is  basic and defines a foliated symplectic form, so it projects to a symplectic form $\hat{\Omega}$ on the transverse manifold $N.$ Assume that the 2-form $\hat{\Omega}$ is the K\"ahler form of the transverse K\"ahler structure $(\hat{g}, \hat{J}).$ Take $D = ker \eta .$  Then it is not difficult to verify that $ (g, \xi, \eta, \phi )$ is a contact metric structure, 
as $\hat{\Omega} (X,Y) = \hat{\Omega} (\hat{J}(X),\hat{J}(Y))$ and $\hat{\Omega} (X,Y) =\hat{g}((X,\hat{J}(Y)).$

\medskip

\noindent
  This equality when lifted to the foliated manifold $(M,{\mathcal F}_{\xi})$ reads

$$ g(X,\phi (Y) ) = d\eta (X,Y) $$

\noindent
The fact that the K\"ahler form $\hat{\Omega}$ is $\hat{J}$ invariant on the foliated manifold $(M,{\mathcal F}_{\xi})$ reads as

$$d \eta (\phi (X), \phi (Y) ) = d\eta (X,Y).$$

\noindent
The condition $d\eta (\phi (X), X ) >0$ for $ 0 \neq X \in D,$ translates itself on the level of the transverse manifold to $\hat{\Omega} ( \hat{J} X, X) > 0$ which follows immediately from the definition of the form $\Omega$ : $\Omega (\hat{J}X,X) = \hat{g} (\hat{J}(X),\hat{J}(X)) > 0$  provided that $X \neq 0.$ 
Therefore strict contact structure $\eta$ whose characteristic foliation is transversely K\"ahler admits a contact metric structure whose 1-form is the contact form $\eta .$ 

\medskip

If the manifold $M$ is compact, the foliated symplectic form $\omega = d\eta$ is basic and defines a non-zero 2-basic cohomology class and the 2n basic form $\omega^n$ defines a non-zero 2n-basic cohomology class, so the characteristic foliation is taut, cf. \cite{Ton}. Therefore we can modify the bundle-like metric $g$ along the tangent bundle to the characteristic foliation to a Riemannian metric $g'$ making the characteristic foliation minimal, i.e. the tangent vector field of unit length in the metric $g'$ is Killing. The modification of the metric did preserve the splitting of the tangent bundle. Therefore on the contact manifold $M$ we have a K-contact structure $(g',\xi ', \eta', \phi)$ whose characteristic foliation is the same ${\mathcal F}_{\xi}.$

\medskip

These considerations can be summed up by the following statement

\medskip
{\it No transverse property can distinguish  K-contact manifolds from  Sasakian manifolds.} Transverse properties can only say that a given foliation is not transversely K\"ahler. However, the characteristic foliation of a  K-contact manifold can be transversely K\"ahler without the structure itself being Sasakian.

\medskip
\noindent
Thus if we want to prove that a given K-contact structure on a compact manifold is not Sasakian (i.e. it is not normal) we have to look for some properties which are not transverse, e.g., it is useless to study properties of the basic cohomology of the characteristic foliation. 

\medskip

The Sasakian version of the Hard Lefschetz Theorem proved by B. Cappelletti Montana et al. provides precisely a true obstruction to "being Sasakian," cf. \cite{Cap}

\begin{theorem}

Let $(M, g, \xi, \eta, \phi )$ be a compact connected  Saasaki manifold of dimension m=2n+1. Then for any  $ 0 \leq p \leq n$ the multiplication by the form $\eta \wedge d\eta^p$ induces an isomorphism between $H^{n-p}(M)$ and  $H^{n+p+1}(M)$.

\end{theorem}

To complement this result the authors constructed two nilmanifolds of dimension 5 and 7, respectively, which are K-contact but do not admit any Sasakian structure. To prove that they use the properties of the cohomology ring which can be drived from the Hard Lefschetz Theorem, cf. \cite{Cap2}.

\medskip
The theorem coupled with these examples demonstrates that the Hard Lefschetz property is an obstruction to being Sasakian for compact manifolds. 

\section{Other cohomology theories}

In search for more cohomological obstructions one can turn to other cohomology theories which have been developed for complex manifolds, for the most recent and up-to-date information see \cite{An}. The foliated versions of several of these cohomology theories have been defined and studied by P. Ra\'zny, cf. 
\cite{Ra}, a Ph.D. student at the Jagiellonian University.

\subsection{Basic Bott-Chern cohomology of foliations}

Let M be a manifold of dimension $m=p+2q,$ endowed with a transversely Hermitian (i.e. transversely holomorphic, posessing a tranverse Hermitian metric) foliation $\mathcal{F}$ of complex codimension q. We can  define the basic de Rham complex (denoted $A^{*}_C(M, \mathcal{F})$) as the subcomplex of the standard de Rham complex of M consisting of basic forms. As in the manifold case the transversly holomorphic structure induces a decomposition of the cotangent spaces into forms of type (0,1) and (1,0), cf. Section 3. The basic Bott-Chern cohomology of $\mathcal{F}$ is defined as

$$
H^{*,*}_{BC}(M, \mathcal{F}):=\frac{Ker(\partial)\cap Ker(\bar{\partial})}{Im(\partial\bar{\partial})}
$$

\noindent
{\bf Remark} 
Complex conjugation induces an antilinear isomorphism:
\begin{equation*}
H^{p,q}_{BC}(M,\mathcal{F})\rightarrow H^{q,p}_{BC}(M,\mathcal{F})
\end{equation*}
in particular:
\begin{equation*}
dim_{\mathbb{C}}H^{p,q}_{BC}(M, \mathcal{F})=dim_{\mathbb{C}}H^{q,p}_{BC}(M, \mathcal{F})
\end{equation*}

\medskip

P. Ra\'zny proves a  decomposition theorem for basic Bott-Chern cohomology using  the operator:
\begin{equation*}
\Delta_{BC}:=(\partial\bar{\partial})(\partial\bar{\partial})^*+(\partial\bar{\partial})^*(\partial\bar{\partial})+
(\bar{\partial}^*\partial)(\bar{\partial}^*\partial)^*+(\bar{\partial}^*\partial)^*(\bar{\partial}^*\partial)
+\bar{\partial}^*\bar{\partial}+\partial^*\partial
\end{equation*}

\noindent
where  $\partial^*$ and $\bar{\partial}^*$ are the adjoint operators  to $\partial$ and $\bar{\partial}$,  respectively, with respect to the Hermitian product, defined by the transverse Hermitian structure, as defined in \cite{E1}. He notices that the operator $\Delta_{BC}$ is transversely elliptic and self-adjoint. To prove ellipticity he uses the fact that the operator projects, on the local quotient manifold, to the manifold version of the $\Delta_{BC}$ operator, which is elliptic, cf. \cite{S1}.

\begin{theorem}(Decomposition of the basic Bott-Chern cohomology)
If M is a compact manifold, endowed with a transversely Hermitian foliation $\mathcal{F}$, then we have the following decomposition:
\begin{equation*}
A^{*,*}(M, \mathcal{F})=Ker(\Delta_{BC})\oplus Im(\partial\bar{\partial})\oplus (Im(\partial^*)+Im(\bar{\partial}^*))
\end{equation*}
In particular:
\begin{equation*}
H^{*,*}_{BC}(M, \mathcal{F})\cong Ker(\Delta_{BC})
\end{equation*}
and the dimension of $H^{*,*}_{BC}(M, \mathcal{F})$ is finite.

\end{theorem}

\subsection{Basic Aeppli cohomology of foliations}

 We define the basic Aeppli cohomology of $\mathcal{F}$ as:
\begin{equation*}
H^{*,*}_{A}(M, \mathcal{F}):=\frac{Ker(\partial\bar{\partial})}{Im(\partial)+Im(\bar{\partial})}
\end{equation*}

\noindent
{\bf Remark} As in the Bott-Chern case complex conjugation induces an antilinear isomorphism:
\begin{equation*}
H^{p,q}_A(M, \mathcal{F})\rightarrow H^{q,p}_A(M, \mathcal{F})
\end{equation*}

To obtain a decomposition theorem for the basic Aeppli cohomology of $\mathcal{F}$ we define a basic  self-adjoint, transversely elliptic  differential operator $\Delta_A$ :

\begin{equation*}
\Delta_A:=\partial\partial^*+\bar{\partial}\bar{\partial}^*+
(\partial\bar{\partial})^*(\partial\bar{\partial})+(\partial\bar{\partial})(\partial\bar{\partial})^*+(\bar{\partial}\partial^*)^*(\bar{\partial}\partial^*)+
(\bar{\partial}\partial^*)(\bar{\partial}\partial^*)^*
\end{equation*}

\noindent
and thus we have 

\begin{theorem}(Decomposition of basic Aeppli cohomology)
Let M be a compact manifold, endowed with a Hermitian foliation $\mathcal{F}$. Then we have the following decomposition
\begin{equation*}
A^{*,*}(M, \mathcal{F})=Ker(\Delta_A)\oplus (Im(\partial)+Im(\bar{\partial}))\oplus Im((\partial\bar{\partial})^*)
\end{equation*}
In particular there is an isomorphism:
\begin{equation*}
H^{*,*}_A(M, \mathcal{F})\cong Ker(\Delta_A)
\end{equation*}
and the dimension of $H^{*,*}_A(M, \mathcal{F})$ is finite.

\end{theorem}

 A duality theorem for basic Bott-Chern and Aeppli cohomology is also true, but we  have to assume that   our foliation is homologically orientable.

\noindent
{\bf Remark} The above condition guaranties, that the following equalities hold for basic $r$-forms:
\begin{equation*}
\partial^*=(-1)^r*\partial*, 
\quad
\bar{\partial}^*=(-1)^r*\bar{\partial}*
\end{equation*}

\begin{corollary}
If M is a compact manifold endowed with a Hermitian, homologicaly orientable foliation $\mathcal{F}$, then the transverse star operator induces an isomorphism:
\begin{equation*}
H^{p,q}_{BC}(M, \mathcal{F})\rightarrow H^{n-p,n-q}_A(M, \mathcal{F})
\end{equation*} 

\end{corollary}

The theorem below is the  main result concerning the basic Bott-Chern and Aeppli cohomologies proved in \cite{Ra}:

\begin{theorem}(Basic Fr\"{o}licher-type inequality)
Let M be a manifold of dimension n, endowed with a transversely holomorphic foliation $\mathcal{F}$ of complex codimension q. Let us assume that the basic Dolbeault cohomology of $\mathcal{F}$ are finitely dimensional. Then, for every $k\in\mathbb{N}$, the following inequality holds:
\begin{equation*}
\sum\limits_{p+q=k}(dim_{\mathbb{C}}(H^{p,q}_{BC}(M, \mathcal{F}))+dim_{\mathbb{C}}(H^{p,q}_A(M, \mathcal{F})))\geq 2dim_{\mathbb{C}}(H^k(M, \mathcal{F},\mathbb{C}))
\end{equation*}
Furthermore, the equality holds for every $k\in\mathbb{N}$, iff $\mathcal{F}$ satisfies the $\partial\bar{\partial}$-lemma (i.e. it's basic Dolbeault double complex satisfies the $\partial\bar{\partial}$-lemma).

\end{theorem}

In the  case when $\mathcal{F}$ is a transversely Hermitian foliation on a closed manifold $M,$  we get the following corollary:

\begin{corollary}
Let F be a transversely Hermitian, homologicaly orientable foliation on a closed manifold M. Then for all $k\in\mathbb{N}$ the following inequality holds:
\begin{equation*}
\sum\limits_{p+q=k}(dim_{\mathbb{C}}(H^{p,q}_{BC}(M, \mathcal{F}))+dim_{\mathbb{C}}(H^{p,q}_A(M, \mathcal{F})))\geq 2dim_{\mathbb{C}}(H^k(M, \mathcal{F},\mathbb{C}))
\end{equation*}
Furthermore, the equality holds for every $k\in\mathbb{N}$, iff $\mathcal{F}$ satisfies the $\partial\bar{\partial}$-lemma (i.e. it's basic Dolbeault double complex satisfies the $\partial\bar{\partial}$-lemma).
\end{corollary}

These properties of new  basic cohomology theories provide new tools to distinguish various transverse structures of Riemannian and  holomorphic foliations. They can be used to check whether a given Riemannian foliation admits a rich transverse holomorphic structure, in partticular whether it is transversely K\"ahler.

\end{document}